# Decentralized Voltage and Power Regulation Control of Excitation and Governor System with Global Asymptotic Stability

Hui Liu, *Member, IEEE*, Junjian Qi, *Member, IEEE*,
Jianhui Wang, *Senior Member, IEEE*, and Peijie Li, *Member, IEEE*,

*Abstract*—The Global Asymptotic Stability (GAS), Voltage Regulation (VR), and Power Regulation (PR) of the excitation and governor control system are of critical importance for power system security. However, simultaneously fulfilling GAS, VR, and PR has not yet been achieved. In order to solve this problem, in this paper, we propose a Lyapunov-based decentralized Control (LBC) for the excitation and governor system of multi-machine power system. A completely controllable linear system is actively constructed to design the time-derivative of the Lyapunov function and GAS is guaranteed by satisfying the condition of GAS in Lyapunov theorem. At the same time, VR and PR are performed by introducing both voltage and power to the feedback. The effectiveness of the proposed method is tested and validated on a six-machine power system.

## NOMENCLATURE

| | |
|---|---|
| $C_{Hi}$ | Power coefficient of HP cylinder, 0.3 |
| $C_{Ii}$ | Power coefficient of IP cylinder, 0.4 |
| $C_{Li}$ | Power coefficient of LP cylinder, 0.3 |
| $E_{fi}$ | Excitation voltage of the excitation system in p.u. |
| $P_{Mi}$ | Mechanical power input in p.u. |
| $P_{ei}$ | Active power of a generator unit in p.u. |
| $P_{0i}$ | Expected value of the active power in p.u. |
| $P_{Hi}$ | Power output of HP cylinder in p.u. |
| $P_{Ri}$ | Power output of re-heater in p.u. |
| $P_{Ii}$ | Power output of IP cylinder in p.u. |
| $P_{Li}$ | Power output of LP cylinder in p.u. |
| $T_{Wi}$ | Water starting time, 1 s |
| $T_{WSi}$ | Time constant of the servomotor of HTG system, 5 s |
| $T_{CSi}$ | Time constant of the servomotor of CG system, 0.2 s |
| $T_{Ci}$ | Time constant of CG system, 0.2 s |
| $T_{Hi}$ | Time constant of HP cylinder, 0.2 s |
| $T_{HSi}$ | Time constant of the servomotor of HP cylinder, 0.2 s |
| $T_{Ri}$ | Time constant of the re-heater, 10 s |
| $T_{Ii}$ | Time constant of IP cylinder, 0.1 s |
| $T_{ISi}$ | Time constant of the servomotor of IP cylinder, 0.2 s |
| $T_{Li}$ | Time constant of LP cylinder, 0.1 s |
| $U_{ti}$ | Generator terminal voltage in p.u. |
| $U_{0i}$ | Expected value of terminal voltage in p.u. |
| $U_{Wi}$ | Opening control signal of the guide vane in p.u. |
| $U_{Ci}$ | Opening control signal of steam valve of CG system in p.u. |
| $U_{Hi}$ | Opening control signal of HP cylinder in p.u. |
| $U_{Ii}$ | Opening control signal of IP cylinder in p.u. |
| $\omega_i$ | Rotor speed in rad/s |
| $\mu_{Wi}$ | Water gate opening in p.u. |
| $\mu_{Ci}$ | Steam valve opening of CG system in p.u. |
| $\mu_{Hi}$ | Steam valve opening of HP cylinder in p.u. |
| $\mu_{Ii}$ | Steam valve opening of IP cylinder in p.u. |

This work was supported in part by National Natural Science Foundation of China (51107054, 51577085) and Senior Professional Science Foundation of Jiangsu University (09JDG008, 10JDG085).

H. Liu is with the School of Electrical and Information Engineering, Jiangsu University, Zhenjiang, China, and is visiting the Energy Systems Division, Argonne National Laboratory, Argonne, IL 60439 USA (e-mail: hughlh@126.com).

J. Qi and J. Wang are with the Energy Systems Division, Argonne National Laboratory, Argonne, IL 60439 USA (e-mails: jqi@anl.gov; jianhui.wang@anl.gov).

P. Li is with the College of Electrical and Engineering, Guangxi University, Nanning, China, and is visiting the Energy Systems Division, Argonne National Laboratory, Argonne, IL 60439 USA (e-mail:lipeijie@gxu.edu.cn)

## I. INTRODUCTION

EXCITATION control and governor control of generating units are important measures for improving transient stability and achieving Voltage Regulation (VR) and Power Regulation (PR) [1]-[16]. For excitation control design, enhancing system stability and maintaining VR are two main goals [3]. In [7], the

Global Asymptotic Stability (GAS) and the VR have been simultaneously addressed by nonlinear control design of the excitation system. However, simultaneously achieving GAS, VR, and PR for both excitation control and governor control is still an open question.

Governor control and excitation control are often considered as two independent controls [8]-[10]. In order to independently discuss excitation control and governor control, assumptions that neglect the mutual influence of the excitation system and governor system must be introduced, which results in lack of the proof on GAS. In [8], a Differential Geometric Control (DGC) method is proposed only for the steam valve governor system to enhance transient stability, where the transient EMF in $q$-axis must be assumed constant. In [9] and [10], decentralized excitation control and steam valve control are discussed by Direct Feedback Linearization (DFL) to improve system stability based on the assumption that the mechanical power and the transient EMF in $q$-axis are constant. Besides, neither VR nor PR can be achieved due to a lack of voltage or power feedback [8]-[10].

By contrast, by considering the mutual interaction between the excitation and governor loops, better transient stability margins can be provided [11]. Taking into account the mutual interaction, DGC is applied to design excitation and governor control [12]-[14]. In [12], a DGC-based decentralized control is proposed for the excitation and governor system of hydraulic generating units. The GAS can be achieved because there is only one equilibrium point for the equivalent linear system. However, since voltage and power are not considered as feedback, VR and PR cannot be performed. In [13], another DGC-based method is developed for large reheat-type turbo-generators. But GAS cannot be achieved because the model of large turbo-generators with High-Pressure (HP), Intermediate-Pressure (IP), and Low-Pressure (LP) cylinders and re-heater dynamics is too complicated to be completely linearized by DGC. In [14], although the DGC-based method achieves GAS and VR for condensing-type turbo-generators, it does not consider PR.

Although GAS, VR, and PR are all very important for power system security, there is no existing research that can achieve all of them at the same time. In this paper, we aim at simultaneously fulfilling GAS, VR, and PR by a Lyapunov-based control (LBC) for the decentralized excitation and governor control system. Our contributions can be summarized as follows:

1) The GAS of a system is guaranteed by satisfying the condition of GAS in Lyapunov theorem based on the design of the eigenvalues of a symmetric real matrix;

2) Voltage and power deviations are introduced to determine the negative definiteness of the time-derivative of Lyapunov function based on a completely controllable linear system that is actively constructed by including voltage and power deviations to perform VR and PR through the control inputs;

3) The proposed LBC method can simultaneously achieve VR, PR, and GAS for different types of generators, such as the hydraulic generators, the condensing-type generators, and the reheat-type generators.

The rest of this paper is organized as follows. In Section II, power system models for different types of generating units are introduced. In Section III, the feedback on the decentralized excitation and governor control is addressed. In Section IV, Lyapunov-based decentralized excitation and governor control is proposed for multi-machine power systems. Simulation results on a six-machine power system are presented in Section V to validate the effectiveness of the proposed control method. Finally, the conclusions are drawn in Section VI.

## II. Models of Excitation and Governor Systems

Generator control systems include the excitation system and governor system. Here, we introduce different types of excitation system and governor system.

### A. Excitation System

Extensive studies have been conducted for the designing of decentralized excitation controllers to enhance power system stability [1]-[7], [17]. In these studies, the excitation control system is usually described by the classical third-order model, for which the mathematical expressions and symbols can be found in many references, such as [1]-[7].

### B. Governor Control System

The Hydraulic Turbine Governor (HTG) system and the Steam Turbine Governor (STG) system are introduced as follows.

**(1) HTG control system**

The HTG control system is used to drive hydro-generator units. It exhibits high-order nonlinear behavior and can be very complex. Without considering the elasticity effect of the water column, the hydraulic turbine of HTG can be described as [12]:

$$\dot{P}_{Mi} = \frac{2}{T_{Wi}}(-P_{Mi} + \mu_{Wi} - T_{Wi}\dot{\mu}_{Wi}). \quad (1)$$

The water-gate servomotor regulating the water gate opening is represented by a first-order inertial system as [12]

$$\dot{\mu}_{Wi} = \frac{1}{T_{WSi}}(-\mu_{Wi} + U_{Wi}). \quad (2)$$

**(2) STG control system**

The STG control system is used to drive turbo-generator units. It usually includes Condensing-type Governor (CG) control system [8] and Reheat-type Governor (RG) control system [16]. The CG system drives small turbo-generators, while the RG system is for large turbo-generators.

a) CG control system

Steam turbine dynamic:

$$\dot{P}_{Mi} = \frac{1}{T_{Ci}}(-P_{Mi} + \mu_{Ci}). \quad (3)$$

The servomotor that regulates the steam flow of the steam turbine can be described by

$$\dot{\mu}_{Ci} = \frac{1}{T_{CSi}}(-\mu_{Ci} + U_{Ci}). \tag{4}$$

b) RG control system

HP cylinder dynamic:
$$\dot{P}_{Hi} = \frac{1}{T_{Hi}}(C_{Hi}\mu_{Hi} - P_{Hi}). \tag{5}$$

The servomotor of the HP cylinder used to regulate the steam flow can be represented by
$$\dot{\mu}_{Hi} = \frac{1}{T_{HSi}}(U_{Hi} - \mu_{Hi}). \tag{6}$$

Re-heater dynamic:
$$\dot{P}_{Ri} = \frac{1}{T_{Ri}}(\frac{P_{Hi}}{C_{Hi}} - P_{Ri}). \tag{7}$$

IP cylinder dynamic:
$$\dot{P}_{Ii} = \frac{1}{T_{Ii}}(C_{Ii}P_{Ri}\mu_{Ii} - P_{Ii}). \tag{8}$$

The servomotor of the IP cylinder is used to regulate the steam valve opening and can be described by
$$\dot{\mu}_{Ii} = \frac{1}{T_{ISi}}(U_{Ii} - \mu_{Ii}). \tag{9}$$

LP cylinder dynamic:
$$\dot{P}_{Li} = \frac{1}{T_{Li}}(\frac{C_{Li}P_{Ii}}{C_{Ii}} - P_{Li}). \tag{10}$$

For the reheat type of generation system, the mechanical power input can be calculated as:
$$P_{Mi} = P_{Hi} + P_{Ii} + P_{Li} \tag{11}$$

where
$$P_{Hi} = C_{Hi}P_{Mi}, \; P_{Ii} = C_{Ii}P_{Mi}, \; P_{Li} = C_{Li}P_{Mi},$$
$$C_{Hi} + C_{Ii} + C_{Li} = 1.$$

## III. FEEDBACK STATEMENT ON EXCITATION AND GOVERNOR CONTROL SYSTEMS

In the control design of the excitation and governor systems, the feedback information is of great importance for achieving control goals. The feedbacks of voltage and power can be used to performing VR and PR [3], while the rotor speed feedback can improve the transient stability of power systems [2]. Besides, some feedback information, such as the water gate opening, should also be included in order to implement PR by regulating the output of active power. In the following, we discuss the feedback information for the excitation system and different types of governor system.

### A. Feedback Information of Excitation Control System

The excitation system is the only way to regulate the generator voltage and its primary objective is to maintain the voltage level at the generator bus. Thus voltage should be used as feedback. Besides, in real power systems the rotor speed is always used as the feedback input of power system stabilizer for excitation control in order to damp system oscillations [18]. Therefore, the feedback information of the excitation system should include
$$\begin{cases} \Delta U_{ti} = U_{ti} - U_{0i} \\ \Delta \omega_i = \omega_i - 1. \end{cases} \tag{12}$$

### B. Feedback Information of Governor Control System

**(1) HTG control system**

The HTG control system consists of the control and actuating equipment used to regulate the power output of hydro- generators by controlling water flow. The key to controlling the water flow is to regulate water-gate opening according to a regulation signal. Therefore, both the power output and the water gate opening are considered as the feedback of nonlinear control as:
$$\begin{cases} \Delta P_{ei} = P_{ei} - P_{0i} \\ \Delta \mu_{Wi} = \mu_{Wi} - P_{ei}. \end{cases} \tag{13}$$

**(2) STG control System**

a) CG control system

For this type of governor system, PR is performed by regulating the steam valve opening. Therefore, the feedback should include the power output and the steam valve opening:
$$\begin{cases} \Delta P_{ei} = P_{ei} - P_{0i} \\ \Delta \mu_{Ci} = \mu_{Ci} - P_{ei}. \end{cases} \tag{14}$$

b) RG control system

For the reheat-type governor control system, PR is controlled by HP, IP, and LP cylinders. From (5), (6), and (11), the objective of the steam valve opening of the HP cylinder can be developed. Similarly, we can deduce the objective of the steam valve opening of the IP cylinder from (8) and (11). Accordingly, we use the following feedback:
$$\begin{cases} \Delta P_{ei} = P_{ei} - P_{0i} \\ \Delta \mu_{Hi} = \mu_{Hi} - P_{ei} \\ \Delta \mu_{Ii} = \mu_{Ii} - P_{ei}/P_{Ri}. \end{cases} \tag{15}$$

## IV. DECENTRALIZED VOLTAGE AND POWER REGULATION CONTROL WITH GAS

Here, GAS, VR, and PR are achieved for the excitation and governor control systems based on Lyapunov theorem.

### A. Lyapunov Function

Based on the above discussion about the three types of governor control systems, we consider a power system with $n^1$ hydro-generators, $n^2$ condensing-type turbo-generators, and $n^3$ reheat-type turbo-generators. By introducing the feedback in (12)-(15), a Lyapunov function is constructed as:

$$V = \frac{1}{2}\sum_{i=1}^{n^1}\left(\Delta\omega_i^2 + \Delta U_{ti}^2 + \Delta P_{ei}^2 + \Delta\mu_{Wi}^2\right)$$
$$+ \frac{1}{2}\sum_{i=1}^{n^2}\left(\Delta\omega_i^2 + \Delta U_{ti}^2 + \Delta P_{ei}^2 + \Delta\mu_{Ci}^2\right) \quad (16)$$
$$+ \frac{1}{2}\sum_{i=1}^{n^3}\left(\Delta\omega_i^2 + \Delta U_{ti}^2 + \Delta P_{ei}^2 + \Delta\mu_{Hi}^2 + \Delta\mu_{Ii}^2\right)$$

The time-derivative of $V$ can be expressed as
$$\dot{V} = \Delta y^T \Delta \dot{y} \quad (17)$$
where
$$\Delta y = [(\Delta y_1^1)^T, \cdots, (\Delta y_{n^1}^1)^T, (\Delta y_1^2)^T, \cdots, (\Delta y_{n^2}^2)^T, (\Delta y_1^3)^T, \cdots, (\Delta y_{n^3}^3)^T]^T$$
and
$$\Delta y_i^1 = [\Delta\omega_i, \Delta U_{ti}, \Delta P_{ei}, \Delta\mu_{Wi}]^T,$$
$$\Delta y_i^2 = [\Delta\omega_i, \Delta U_{ti}, \Delta P_{ei}, \Delta\mu_{Ci}]^T,$$
$$\Delta y_i^3 = [\Delta\omega_i, \Delta U_{ti}, \Delta P_{ei}, \Delta\mu_{Hi}, \Delta\mu_{Ii}]^T.$$

As in (17), the negative definiteness of $\dot{V}$ depends on the differential trajectory $\Delta\dot{y}$. Therefore, we should design the differential trajectory $\Delta\dot{y}$ through the control inputs to guarantee that $\dot{V}$ is negative definite.

### B. Design of Differential Trajectory $\Delta\dot{y}$

In order to design the differential trajectory $\Delta\dot{y}$ through control inputs, we should first deduce the relationship between $\Delta\dot{y}$ and the control inputs.

The time-derivative of the terminal voltage in $\Delta y$ of (17) can be deduced for each generating unit as:
$$\Delta\dot{U}_{ti} = c_i^E + d_i^E E_{fi} \quad (18)$$
where $c_i^E$ and $d_i^E$ describe the relationship between $\Delta\dot{U}_{ti}$ and the excitation voltage $E_{fi}$, and more details can be found in [7].

For the HTG control system, the time-derivative of $\Delta\mu_{Wi}$ in $\Delta y$ can be obtained from (2) and (13) as
$$\Delta\dot{\mu}_{Wi} = c_i^W + d_i^W U_{Wi} \quad (19)$$
where
$$c_i^W = -\frac{\mu_{Wi}}{T_{WSi}} - \dot{P}_{ei}, \quad d_i^W = \frac{1}{T_{WSi}}.$$

For the condensing-type governor control system, with (4) and (14), the time-derivative of $\Delta\mu_{Ci}$ in $\Delta y$ is:
$$\Delta\dot{\mu}_{Ci} = c_i^C + d_i^C U_{Ci} \quad (20)$$
where
$$c_i^C = -\frac{\mu_{Ci}}{T_{CSi}} - \dot{P}_{ei}, \quad d_i^C = \frac{1}{T_{CSi}}.$$

Likewise, $\Delta\dot{\mu}_{Hi}$ and $\Delta\dot{\mu}_{Ii}$ can be derived from (6), (9), and (15) for a reheat-type governor system as:
$$\begin{cases} \Delta\dot{\mu}_{Hi} = c_i^H + d_i^H U_{Hi} \\ \Delta\dot{\mu}_{Ii} = c_i^I + d_i^I U_{Ii} \end{cases} \quad (21)$$
where

$$c_i^H = -\frac{\mu_{Hi}}{T_{HSi}} - \dot{P}_{ei}, \quad d_i^H = \frac{1}{T_{HSi}}, \quad c_i^I = -\frac{\mu_{Hi}}{T_{HSi}} - \dot{P}_i, \quad d_i^I = \frac{1}{T_{ISi}}$$
and
$$P_i = \frac{P_{ei}}{P_{Ri}}.$$

With (18)-(21), the equations in (12)-(15) are used to actively construct a completely controllable linear system, as follows:
$$\Delta\dot{y} = A\Delta y + B(c + Du) \quad (22)$$
where
$$A = diag\left(A_1^1, \cdots, A_{n^1}^1, A_1^2, \cdots, A_{n^2}^2, A_1^3, \cdots, A_{n^3}^3\right)$$
$$B = diag\left(B_1^1, \cdots, B_{n^1}^1, B_1^2, \cdots, B_{n^2}^2, B_1^3, \cdots, B_{n^3}^3\right)$$
$$c = \left[(c_1^1)^T, \cdots, (c_{n^1}^1)^T, (c_1^2)^T, \cdots, (c_{n^2}^2)^T, (c_1^3)^T, \cdots, (c_{n^3}^3)^T\right]^T$$
$$u = \left[(u_1^1)^T, \cdots, (u_{n^1}^1)^T, (u_1^2)^T, \cdots, (u_{n^2}^2)^T, (u_1^3)^T, \cdots, (u_{n^3}^3)^T\right]^T$$
$$D = diag\left(d_1^1, \cdots, d_{n^1}^1, d_1^2, \cdots, d_{n^2}^2, d_1^3, \cdots, d_{n^3}^3\right)$$
and
$$A_i^1 = \begin{bmatrix} a_{i,1}^1 & a_{i,2}^1 & 0 & 0 \\ 0 & 0 & 0 & 0 \\ 0 & 0 & a_{i,3}^1 & a_{i,4}^1 \\ 0 & 0 & 0 & 0 \end{bmatrix}, \quad B_i^1 = \begin{bmatrix} 0 & 0 \\ 1 & 0 \\ 0 & 0 \\ 0 & 1 \end{bmatrix},$$
$$A_i^2 = \begin{bmatrix} a_{i,1}^2 & a_{i,2}^2 & 0 & 0 \\ 0 & 0 & 0 & 0 \\ 0 & 0 & a_{i,3}^2 & a_{i,4}^2 \\ 0 & 0 & 0 & 0 \end{bmatrix}, \quad B_i^2 = \begin{bmatrix} 0 & 0 \\ 1 & 0 \\ 0 & 0 \\ 0 & 1 \end{bmatrix},$$
$$A_i^3 = \begin{bmatrix} a_{i,1}^3 & a_{i,2}^3 & 0 & 0 & 0 \\ 0 & 0 & 0 & 0 & 0 \\ 0 & 0 & a_{i,3}^3 & a_{i,4}^3 & 0 \\ 0 & 0 & 0 & 0 & 0 \\ 0 & 0 & 0 & 0 & 0 \end{bmatrix}, \quad B_i^3 = \begin{bmatrix} 0 & 0 & 0 \\ 1 & 0 & 0 \\ 0 & 0 & 0 \\ 0 & 1 & 0 \\ 0 & 0 & 1 \end{bmatrix}.$$
$$u_i^1 = [E_{fi}, U_{Wi}]^T, \quad u_i^2 = [E_{fi}, U_{Ci}]^T, \quad u_i^3 = [E_{fi}, U_{Hi}, U_{Ii}]^T$$
$$c_i^1 = [c_i^E, c_i^W]^T, \quad c_i^2 = [c_i^E, c_i^C]^T, \quad c_i^3 = [c_i^E, c_i^H, c_i^I]^T$$
$$d_i^1 = \begin{bmatrix} d_i^E & 0 \\ 0 & d_i^W \end{bmatrix}, \quad d_i^2 = \begin{bmatrix} d_i^E & 0 \\ 0 & d_i^C \end{bmatrix}, \quad d_i^3 = \begin{bmatrix} d_i^E & 0 & 0 \\ 0 & d_i^H & 0 \\ 0 & 0 & d_i^I \end{bmatrix}.$$

Let $v = c + Du$, (22) can be rewritten as
$$\Delta\dot{y} = A\Delta y + Bv \quad (23)$$

***Remark 1:*** *Equation (23) is designed to regulate the differential trajectory $\Delta\dot{y}$ through the control inputs. We can actively construct such a system, because control can be viewed a special force presented by the designers to achieve the goals of the control design. According to linear control theory, the system illustrated in (23) can be a completely controllable linear system by the constants $a_{i,j}^1$, $a_{i,j}^2$ and $a_{i,j}^3$ ($j = 1,2,3,4$).*

***Remark 2:*** *With complete controllability, the poles of (23) can be arranged arbitrarily through virtual inputs, by applying*

linear control method, which means that the trajectory of time-derivative $\Delta \dot{y}$ can be controlled by virtual inputs based on pole arrangements. Therefore, we can use (23) to regulate the negative definiteness of (17).

***Remark 3***: It is noted that (23) is constructed with the feedback information such as voltage and power deviations and thus we introduce both voltage and power deviations to determine the negative definition of the time-derivative of *Lyapunov* function. Therefore, GAS is closely related to VR and PR, and GAS, PR, and VR are considered simultaneously.

For (23), the feedback can be solved as
$$v = -K \Delta y \qquad (24)$$
where
$$K = diag(k_1^1, \cdots, k_{n^1}^1, k_1^2, \cdots, k_{n^2}^2, k_1^3, \cdots, k_{n^3}^3)$$
and
$$k_i^1 = \begin{bmatrix} k_{i,1}^1 & k_{i,2}^1 & 0 & 0 \\ 0 & 0 & k_{i,3}^1 & k_{i,4}^1 \end{bmatrix}$$
$$k_i^2 = \begin{bmatrix} k_{i,1}^2 & k_{i,2}^2 & 0 & 0 \\ 0 & 0 & k_{i,3}^2 & k_{i,4}^2 \end{bmatrix}$$
$$k_i^3 = \begin{bmatrix} k_{i,1}^3 & k_{i,2}^3 & 0 & 0 & 0 \\ 0 & 0 & k_{i,3}^3 & k_{i,4}^3 & 0 \\ 0 & 0 & 0 & 0 & k_{i,5}^3 \end{bmatrix}.$$

The decentralized excitation and governor control can be obtained from $v = c + Du$ and (24) as
$$u = D^{-1}(-K\Delta y - c) \qquad (25)$$
which is
$$u_i^1 = \begin{bmatrix} E_{fi} \\ U_{Wi} \end{bmatrix} = \begin{bmatrix} (-k_{i,1}^1 \Delta \omega_i - k_{i,2}^1 \Delta U_{ti} - c_i^E)/d_i^E \\ (-k_{i,3}^1 \Delta P_{ei} - k_{i,4}^1 \Delta \mu_{Wi} - c_i^W)/d_i^W \end{bmatrix}$$
$$u_i^2 = \begin{bmatrix} E_{fi} \\ U_{Ci} \end{bmatrix} = \begin{bmatrix} (-k_{i,1}^2 \Delta \omega_i - k_{i,2}^2 \Delta U_{ti} - c_i^E)/d_i^E \\ (-k_{i,3}^2 \Delta P_{ei} - k_{i,4}^2 \Delta \mu_{Ci} - c_i^C)/d_i^C \end{bmatrix}$$
$$u_i^3 = \begin{bmatrix} E_{fi} \\ U_{Hi} \\ U_{Ii} \end{bmatrix} = \begin{bmatrix} (-k_{i,1}^3 \Delta \omega_i - k_{i,2}^3 \Delta U_{ti} - c_i^E)/d_i^E \\ (-k_{i,3}^3 \Delta P_{ei} - k_{i,4}^3 \Delta \mu_{Hi} - c_i^H)/d_i^H \\ (-k_{i,5}^3 \Delta \mu_{Ii} - c_i^I)/d_i^I \end{bmatrix}.$$

### C. Justification on GAS

Substituting (22) and (25) into (17), the time-derivative of the *Lyapunov* function can be rewritten as:
$$\dot{V} = \Delta y^T \Phi \Delta y \qquad (26)$$
where
$$\Phi = A - BK. \qquad (27)$$

We define a symmetric real matrix $\Psi = \Phi + \Phi^T$ whose negative definiteness is equivalent to that of $\Phi$ according to matrix theory. By properly choosing the coefficient $K$, we can make sure that all of the eigenvalues of $\Psi$ are negative real numbers and thus $\Psi$ is negative definite. Then for any $\Delta y \neq 0$, there is:

$$\dot{V} = \Delta y^T \Phi \Delta y < 0. \qquad (28)$$

Denote the state deviation vector of the dynamic equations of a power system by $\Delta x$. If $\Delta x \neq 0 \Rightarrow \Delta y \neq 0$, we can get
$$\dot{V} = \Delta y^T \Phi \Delta y < 0 \text{ for any } \Delta x \neq 0. \qquad (29)$$

With (29), GAS can be accomplished according to the *Lyapunov* theorem. Therefore, we should prove that for any $\Delta x \neq 0$, there is $\Delta y \neq 0$.

For the sake of convenience, we arrange $\Delta x$ and $\Delta y$ as:
$$\begin{cases} \Delta x = [\Delta \delta_1, \Delta E'_{q1}, \Delta \omega_1, \cdots, \Delta \delta_n, \Delta E'_{qn}, \Delta \omega_n, \Delta x_g]^T \\ \Delta y = [\Delta U_{t1}, \Delta P_{e1}, \Delta \omega_1, \cdots, \Delta U_{tn}, \Delta P_{en}, \Delta \omega_n, \Delta y_g]^T \end{cases} \qquad (30)$$
where $n = n^1 + n^2 + n^3$, and $\Delta x_g$ and $\Delta y_g$ are the state vector and feedback vector of the governor systems.

Then we need to show that for any element of $\Delta x$ ($\Delta \delta_i$, $\Delta E'_{qi}$, $\Delta \omega_i$, or $\Delta x_g$) not equal to zero, there is $\Delta y \neq 0$. Specifically, 1) We can easily get $\Delta \omega_i \neq 0 \Rightarrow \Delta y \neq 0$ since $\Delta \omega_i$ is also an element of $\Delta y$; 2) Considering the models and physical characteristics of the governor systems, we can deduce $\Delta x_g \neq 0 \Rightarrow \Delta y_g \neq 0$. For example, when $\Delta \mu_{Ci}$ in $\Delta x_g$ is not equal to zero, $\Delta P_{ei}$ in $\Delta y_g$ will not be equal to zero, that is, the change of the steam valve opening will result in the change of power output; 3) When $\Delta \delta_i$ or $\Delta E'_{qi}$ is not equal to zero, $\Delta U_{ti}$ and $\Delta P_{ei}$ cannot be equal to zero at the same time because in that case the generators will not be controllable.

Therefore, for any $\Delta x \neq 0$, there is $\Delta y \neq 0$. Thus (29) holds and GAS is guaranteed.

### D. Statement on Performing GAS, VR, and PR

The *Lyapunov* function in (16) is constructed by using a quadratic form of the feedback in (12)-(15). The time-derivative of the *Lyapunov* function is also designed as a quadratic form of the above feedback through control inputs based on the design of a differential trajectory. Therefore, GAS is closely related to the tracking errors of the system variables such as terminal voltage and active power. While achieving the stability of a system, the tracking deviations of voltage and active power are also decreasing. When the system finally stabilizes in a steady state due to the feedback control, the tracking errors of voltage and power will be also eliminated. In such a manner, GAS, VR, and PR are achieved simultaneously.

## V. SIMULATION RESULTS

Here, we present results to validate the effectiveness of the proposed control approach. All simulations are performed in Matlab.

### A. System Description

To demonstrate a nonlinear control design, a power system with 3 generators is often used, as in [2], [9], [10], and [17]. However, in order to cover more types of generators, here we consider a six-machine power system illustrated in Fig. 1, for which Generator 1 represents the infinite bus, Generator 6 is a synchronous condenser, Generators 2 and 5 are large reheat-type generators, Generator 3 is a hydro-generator, and Generator 4 is a small condensing-type generator. More details about the parameters of this system can be found in [8].

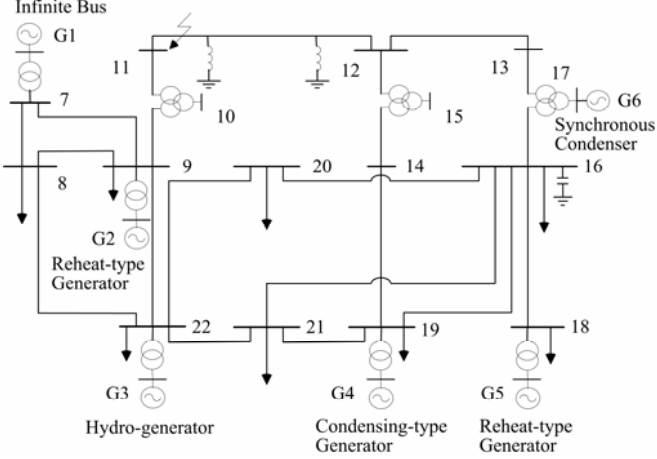

Fig. 1 A six-machine power system

The physical limits of excitation voltages for the excitation control system are set as

$$-5 \leq E_{f_i} \leq 5 \quad (i=2,\cdots,5).$$

The physical limits of the governor control system are:
$0 \leq \mu_{Hi} \leq 7$, $0 \leq \mu_{Ii} \leq 1.1$, $0 \leq \mu_{Ci} \leq 0.8$, $0 \leq \mu_{Wi} \leq 6$.

### B. Parameter Calculation

For generators 2-5, the models in (23) have a total of eighteen orders. By using the parameters in Table I, we have

$$\text{rank}([\boldsymbol{B} \mid \boldsymbol{AB} \mid \boldsymbol{A}^2\boldsymbol{B} \mid \cdots \mid \boldsymbol{A}^{17}\boldsymbol{B}]) = 18,$$

and thus the linear system in (23) is completely controllable

When one considers the feedback gains in Table II, all of the eigenvalues of the matrix $\boldsymbol{\Psi}$ in Section IV.C are negative real number and thus $\boldsymbol{\Psi}$ is negative definite. Consequently, the matrix $\boldsymbol{\Phi}$ is also negative definite. Therefore, GAS can be guaranteed because the condition in (28) and (29) is satisfied.

As in Table II, we use the same gains for all of the generators. System performance may be improved by choosing different feedback gains for different generators, which, however, is out of the scope of this paper.

TABLE I
PARAMETERS OF MATRIX $\boldsymbol{A}$ IN (22)

| $a_{i,1}^1, a_{i,1}^2, a_{i,1}^3$ | $a_{i,2}^1, a_{i,2}^2, a_{i,2}^3$ | $a_{i,3}^1, a_{i,3}^2, a_{i,3}^3$ | $a_{i,4}^1, a_{i,4}^2, a_{i,4}^3$ |
|---|---|---|---|
| -300 | -300 | -10 | 10 |

TABLE II
GAINS OF NONLINEAR FEEDBACKS IN (25)

| $k_{i,1}^1, k_{i,1}^2, k_{i,1}^3$ | $k_{i,2}^1, k_{i,2}^2, k_{i,2}^3$ | $k_{i,3}^1, k_{i,3}^2, k_{i,3}^3$ | $k_{i,4}^1, k_{i,4}^2, k_{i,4}^3$ | $k_{i,5}^3$ |
|---|---|---|---|---|
| -400 | 30 | 5 | 5 | 5 |

### C. Simulation Analysis

In order to demonstrate the effectiveness of the proposed Lyapunov-based Control (LBC), we compare it with the DGC-based methods [12]-[14]. Note that LBC can be used for all types of generators while each DGC-based method can only be used for one specific type of generator. The DGC-based method in [13] for large reheat-type generators is called DGC-R for convenience. The DGC-based method in [12] for hydro-generators is called DGC-H. The DGC-based method in [14] for small condensing-type generators is called DGC-C. Therefore, generators 2-5 are assigned DGC-R, DGC-H, DGC-C, and DGC-R, respectively.

Although DGC-R, DGC-H, and DGC-C are all proposed based on DGC, their control designs are significantly different. In Table III, we show whether or not each of them can achieve VR, PR, or GAS. As discussed in Section IV, the LBC, however, can simultaneously achieve VR, PR, and GAS for all types of generators.

TABLE III
COMPARISON OF DGC-BASED METHODS

|  | Gen. 2 | Gen. 3 | Gen. 4 | Gen. 5 |
|---|---|---|---|---|
| Method | DGC-R | DGC-H | DGC-C | DGC-R |
| VR | Yes | No | Yes | Yes |
| PR | Yes | No | No | Yes |
| GAS | No | Yes | Yes | No |

**(1) Performing VR**

At 0.25 s, there is a step change of voltage references from the initial condition to the regulation target, as shown in Table IV. By LBC and DGC-R, the expected voltage can be achieved, as illustrated in Fig. 2(a). From Fig. 2(b), power outputs of generators can also be maintained at initial values by LBC and DGC-R.

TABLE IV
SIMULATION SCENARIOS FOR PERFORMING VR

|  | Terminal Voltage (p.u.) | | | |
|---|---|---|---|---|
|  | Gen. 2 | Gen. 3 | Gen. 4 | Gen. 5 |
| Initial condition | 1.00 | 1.00 | 1.00 | 1.00 |
| Regulation target | 1.025 | 0.975 | 1.025 | 0.975 |

Although the power output of generator 3 almost keeps constant by DGC-H, as in Fig. 2 (b), its terminal voltage deviates from the initial value due to the lack of voltage feedback and the change of power flow, as shown in Fig. 2 (a). From Fig. 2, by DGC-C, the terminal voltage of generator 4 can be regulated to the expected value, but the active power output deviates from the initial value due to a lack of power feedback.

**(2) Performing PR**

At 0.25 s, there is a step change of power references, as in Table V. The dynamic response is shown in Fig. 3.

By LBC and DGC-R, PR can also be achieved while maintaining the initial voltage. By DGC-H, PR cannot be performed, and terminal voltage also deviates from the initial value, because of it does not have voltage or power feedback. For DGC-C, although the terminal voltage can be kept at initial value, the expected output of active power cannot be achieved because of a lack of power feedback, as shown in Fig. 3(b).

TABLE V
SIMULATION SCENARIOS FOR PERFORMING PR

| | Active Power (p.u.) | | | |
|---|---|---|---|---|
| | Gen. 2 | Gen. 3 | Gen. 4 | Gen. 5 |
| Initial condition | 6.00 | 3.10 | 0.60 | 4.30 |
| Regulation target | 5.00 | 4.00 | 0.80 | 5.00 |

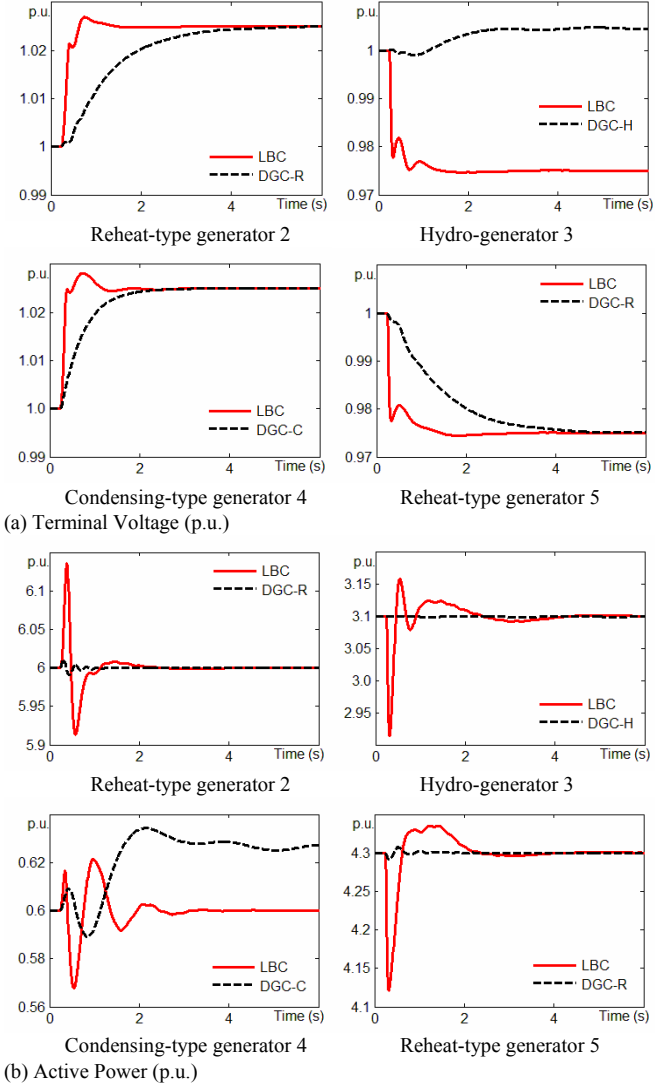

(a) Terminal Voltage (p.u.)

(b) Active Power (p.u.)
Fig. 2 Dynamic response for performing VR

**(3) Three-phase short circuit fault**

At 0.5 s, a three-phase short circuit fault is applied at the beginning of line 11-12 in Fig. 1, and the line is removed after 0.15 second. As shown in Fig. 4(a), except for the voltage deviation of generator 3 resulting from the lack of voltage feedback in DGC-H, the terminal voltages of the other generators can be maintained by using LBC and DGC-based methods.

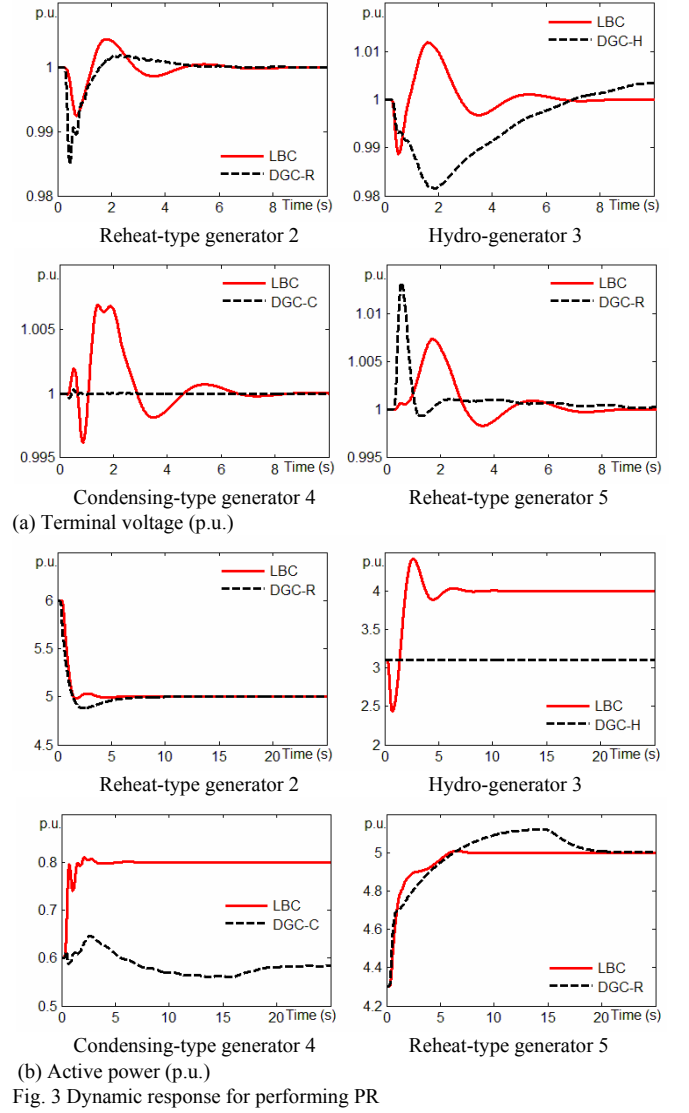

(a) Terminal voltage (p.u.)

(b) Active power (p.u.)
Fig. 3 Dynamic response for performing PR

As in Fig. 4(b), LBC is much more effective than the DGC-based methods in damping system oscillations caused by the three-phase short circuit fault.

In order to further explore the benefits of the proposed controller in enhancing system stability, Critical Clear Time (CCT) is calculated by trial and error by considering different line faults and is listed in Table VI. When a three-phase short circuit fault occurs at line 11-12, the CCT of LBC is the same as that of the DGC-based methods. However, for the other faults, the LBC is much more superior to the DGC-based methods in improving system stability. This is because DGC-R cannot achieve GAS. With DGC-R, the system stability can be enhanced in some cases when appropriate feedback gains are used, but this does not hold for the other cases because GAS cannot be theoretically guaranteed.

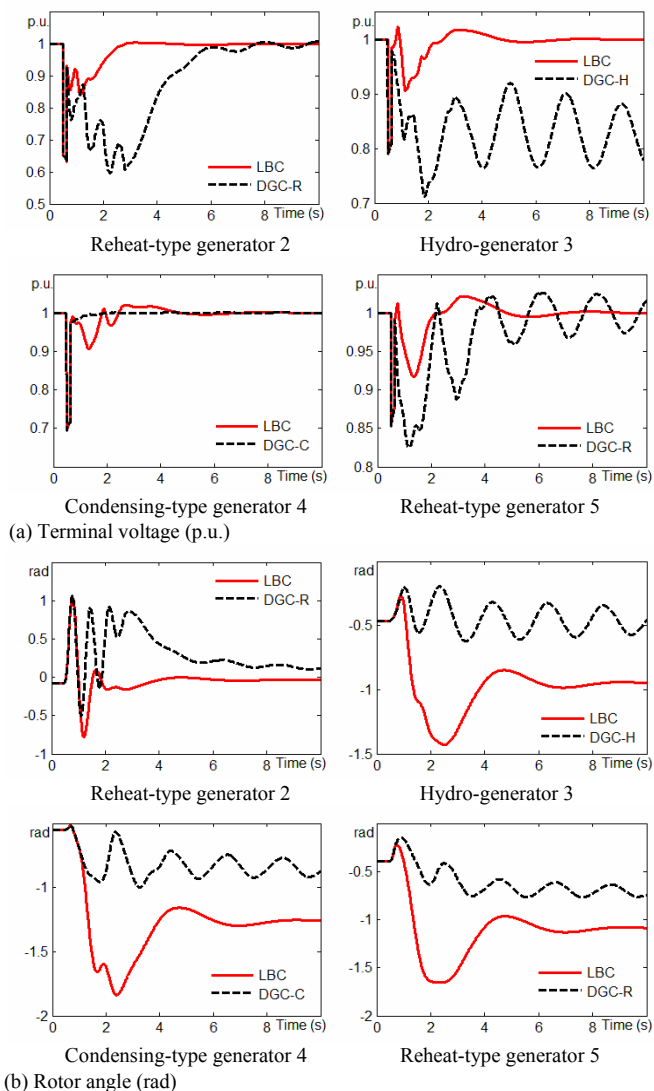

(a) Terminal voltage (p.u.)

(b) Rotor angle (rad)

Fig. 4 Dynamic response to a three phase fault

TABLE VI
CCT COMPARISONS FOR EXCITATION CONTROL AND GOVERNOR CONTROL METHODS

| Fault Bus | Line Removed | Critical Clearing Time (CCT) | |
|---|---|---|---|
| | | LBC (s) | DGC-based methods (s) |
| Bus 11 | Line 11-12 | 0.25 | 0.25 |
| Bus 12 | Line 12-13 | 0.54 | 0.38 |
| Bus 9 | Line 9-22 | 0.20 | 0.16 |

## VI. CONCLUSION

VR, PR, and GAS are important for power system security and thus need to be considered in the excitation and governor control design. The GAS of a power system can be achieved by using some advanced control methods (such as DGC) [12], if VR and PR are not considered. However, to the best of our knowledge, simultaneously fulfilling GAS, VR and PR in the excitation and governor control design has not been achieved.

In this paper, we propose a decentralized excitation and governor controller to tackle this challenge. The GAS of the power system is achieved by the proposed Lyapunov-based controller, and at the same time, voltage and power regulation are performed by introducing both voltage and power into the nonlinear feedback. Simulation results on a six-machine power system demonstrate and validate the effectiveness of the proposed control method.